\documentclass[onenumber,plainnumbering,english,11pt,standalone]{amsart}
\usepackage[utf8]{inputenc}
\usepackage[square,numbers]{natbib}
\usepackage[mathscr]{eucal}
\usepackage{amssymb,amsmath,amscd,mathrsfs,paralist,amsthm}
\usepackage{mathrsfs}
\usepackage{bm}
\usepackage[T1]{fontenc}
\usepackage{tikz}
\usetikzlibrary{shapes.geometric, arrows}                      
\usepackage{tikz}
\usepackage{bm}
\usepackage{graphicx}
\usepackage{pgfplots}
\pgfplotsset{compat=1.15}
\usetikzlibrary{positioning, shapes.geometric}
\usepgfplotslibrary{colormaps}
 
\usepackage{tikz-cd}
\usetikzlibrary{arrows, matrix}
\usepackage[all,cmtip]{xy}
\usepackage[colorlinks,linkcolor=black,urlcolor=black]{hyperref}
\hypersetup{
     colorlinks   = true,
     citecolor    = black}
\usepackage{pgfplots}
\pgfplotsset{compat=newest}
\usepackage[font=small,labelsep=none]{caption}
\usetikzlibrary{matrix}
\usetikzlibrary{cd}
\usepackage[colorlinks,linkcolor=black,urlcolor=black]{hyperref}
\hypersetup{
     colorlinks   = true,
     citecolor    = black}
\usepackage[all,cmtip]{xy}
\theoremstyle{plain}
\newtheorem{theorem}{Theorem}[section]

\newtheorem{lemma}[theorem]{Lemma}
\newtheorem{proposition}[theorem]{Proposition}
\newtheorem{corollary}[theorem]{Corollary}
\theoremstyle{definition}
\newtheorem{definition}[theorem]{Definition}

\theoremstyle{remark}

\makeatletter
\let\amstexbig\big
\def\newbig#1{%
  \ifx#1|%
    \expandafter\@firstoftwo
  \else
    \expandafter\@secondoftwo
  \fi
  {\big@bar}%
  {\amstexbig{#1}}%
}
\AtBeginDocument{\let\big\newbig}
\def\big@bar{\bBigg@{1.1}|}
\makeatother

\DeclareMathOperator{\image}{Im}
\DeclareMathOperator{\Id}{Id}

\DeclareMathOperator{\coker}{Coker}
\DeclareMathOperator{\rk}{rk}
\DeclareMathOperator{\sol}{Sol}
\DeclareMathOperator{\FRAC}{Frac}
\numberwithin{equation}{section}
\newcommand{\C}{\mathbb C}

\newcommand{\theoref}[1]{Theorem~\ref{#1}}
\newcommand{\propref}[1]{Proposition~\ref{#1}}

\newcommand{\lemref}[1]{Lemma~\ref{#1}}

\newcommand{\secref}[1]{Section~\ref{#1}}
\begin{document}
\title{Holomorphic projective connections on surfaces}
\author{Oumar Wone}
\address{Oumar Wone}
    \email{wone@chapman.edu}   
\begin{abstract}
We study complex analytic (possibly singular) projective connections on the plane. We characterize some of them in terms of their families of integral curves. We also give a beginning of classification of second order odes polynomial in the first and second derivatives, and with holomorphic coefficients. \end{abstract}
 \maketitle
 \tableofcontents
\section{Introduction}
\label{chapab}
The differential equation describing the unparameterized geodesics of a projective connection on a surface appeared in mathematics since at least the introduction of the notion of a geodesics of an affine connection on a surface. Indeed it is ubiquitous in classical differential geometry and has been studied by \cite{cartan1924, liouville1889, liouville18892, tresse1894, lie1883}, among others; for a modern presentation see \cite{olver}. Recently \cite{bryant2009, krug2008, luza} there has been a renewal of interest in the study of projective connections on surfaces since its interpretation in terms of twistors by \cite{hitchin1982, lebrun1, kamran}. The interest in projective connections stems from the fact that the differential equation of the geodesics of a (holomorphic) metric on a (holomorphic) surface is described by a projective connection, i.e. by an ordinary differential equation of the form
\begin{equation}
\label{abelseecond}
u^{\prime\prime}+a_1(x,u)u^{\prime3}+3a_2(x,u)u^{\prime2}+3a_3(x,u)u^\prime+a_4(x,u)=0,\, u^\prime=\frac{du}{dx},\,u^{\prime\prime}=\frac{d^2u}{dx^2}
\end{equation}
where $a_i(x,u)\in\C\{x,u\}$ for $1\leqslant i\leqslant4$, and also from the fact that the equivalence classes (under point transformations) of the projective connections is related to the problem of local embedding of a rational curve $\C P_1$ in a holomorphic surface, as a curve with normal bundle of degree $1$.

In this paper we study holomorphic projective connections, and to a lesser extent, second order ordinary differential equations. More precisely we, in \secref{chapab22}, firstly associate, using the Spencer-Goldschmidt theory of overdetermined system of partial differential equations, to a well-chosen system of three second order partial differential equations in two independent variables, a vector bundle of rank three over $(\C^2,0)$ with a non-necessarily flat connection, and such that the flatness of the connection is equivalent to the complete integrability of that system (that is it satisfies all its integrability conditions). As a byproduct we deduce from this some consequences about flat projective connections, that is projective connections which can be linearized to $u^{\prime\prime}=0$, under point transformations. Then we give in \secref{scheffers} a characterization of some meromorphic projective connections in terms of the geometry of a pencil of conics. Finally in \secref{kasner} we study the differential geometry of the tuple $(u^\prime,u^{\prime\prime})$ of coordinates of the second order jet space $J_2(\mathcal O)$ and deduce from this some necessary condition of equivalence, under point transformations, of second order differential equations polynomial in $u^\prime$, $u^{\prime\prime}$, including results on the differential equations defining projective connections. For a different point of view, concerning \secref{chapab22}, see \cite{eastwood2020}.
\section{Holomorphic projective connections and overdetermined systems}
\label{chapab22}
\begin{theorem}
\label{sec1}
Consider the system of partial differential equations
\begin{equation}
\label{abel40}
\begin{split}
&z^{(2,0)}+p_3z^{(1,0)}+q_3z^{(0,1)}+r_3z=0\\
&z^{(1,1)}+p_2z^{(1,0)}+q_2z^{(0,1)}+r_2z=0\\
&z^{(0,2)}+p_1z^{(1,0)}+q_1z^{(0,1)}+r_1z=0\\
&(p_i)_{1\leqslant i\leqslant 3}, (q_i)_{1\leqslant i\leqslant 3}, (r_i)_{1\leqslant i\leqslant 3}\in\C\{x,y\},\,\,z^{(l,j)}=\dfrac{\partial^{l+j}z}{\partial x^l\partial u^j}. \end{split}
\end{equation}
Then we can associate to it a vector bundle of rank three with connection on $(\C^2,O)$, where $O\in\C^2$ is a point that is fixed. This connection is flat if and only if the differential system \eqref{abel40} satisfies all its six integrability conditions. Moreover when the connection is flat its sheaf of horizontal sections is a local system of rank three isomorphic to the solution space of the partial differential equations \eqref{abel40}.
\end{theorem}
\begin{proof}
To give the proof we need some preliminary notions from the Spencer-Goldschmidt theory of overdetermined systems of partial differential equations \cite{bryant1991, quillen1964, spencer1969, henaut2004, seiler2010}, which we adapt to our situation. Consider the system of partial differential equations \eqref{abel40}. Its left hand side is a linear partial differential operator of order $2$ on $\mathcal{O}:=\C\{x,y\}$ with values in $\mathcal{O}^3=(\C\{x,,y\})^3$. Namely we have a map
$$\theta:\mathcal{O}\longrightarrow\mathcal{O}^3.$$
Its corresponding morphism of $\mathcal{O}$-modules 
$$P_0:J_2(\mathcal{O})\longrightarrow\mathcal{O}^3$$
satisfies $P_0\circ j_2=\theta$ where for $z\in\mathcal{O}$ we have $$j_2(z)=(z,z_x,z_u,z_{xx},z_{xu},z_{uu})=:(z,p,q,r,s,t).$$ $j_2$ is from $\mathcal{O}$ to $J_2(\mathcal{O})$ and $(z,z_x,z_u,z_{xx},z_{xu},z_{uu})$ provide coordinates on $J_2(\mathcal{O})$, which is called the space of second order jet of functions. Furthermore $p_0$ is the surjective map onto $\mathcal{O}^3$ given by
\begin{equation}
\label{abel400}
\begin{split}
P_0\colon&J_2(\mathcal{O})\longrightarrow\mathcal{O}^3\\
&(z,z_x,z_u,z_{xx},z_{xu},z_{uu})\mapsto \left(\begin{array}{c} z_{xx}+p_3z_{x}+q_3z_{u}+r_3z\\ z_{xu}+p_2z_{x}+q_2z_{u}+r_2z\\ z_{uu}+p_1z_{x}+q_1z_{u}+r_1z \end{array}\right).\end{split}
\end{equation}
We use the classical notations as in \cite{spencer1969}, \cite{bryant1991}, \cite{quillen1964}, with some slight modifications. Let $P_k:J_{k+2}(\mathcal{O})\to J_k(\mathcal{O}^3)$ for $k\geqslant0$, be the $k$-th prolongation of $P_0$ obtained from it by successive derivations (see \cite[Def.1.2.4]{spencer1969}), and $R_{k}$ the kernel of $p_k$. Besides we have the exact sequence on symbols \cite[p.~185]{spencer1969}
\begin{equation}
\label{abel402}
\begin{tikzcd}
  0 \arrow[r]& g_k \arrow[r]& S_{k+2}(\mathcal{O}) \arrow[r, "\sigma_k" ] & S_k(\mathcal{O}^3)  \arrow[r, "\tau_{k-1}" ] &  \coker(\sigma_k)\rar&0 \end{tikzcd}.
  \end{equation}
 Here $S_{k+1}(\mathcal O)$ designates the kernel of the natural morphism from $J_{k+1}(\mathcal O)\to J_k(\mathcal O)$ and we have a similar definition for $S_{k+1}(\mathcal O^3)$. 
  
  With these notations and using among others the snake lemma (for the existence of the map $\nu_k$ below) \cite[p.~27-28]{quillen1964} we get the exact commutative sequence of $\mathcal{O}$-modules
  \begin{equation}
  \label{abel403}
\begin{tikzcd}
  &   0\ar[d] & 0\ar[d]&0\ar[d]\\
  0 \arrow[r] & g_{k+1} \arrow[d, ""] \arrow[r, ""] & S_{k+3}(\mathcal{O}) \arrow[d, ""] \arrow[r, "\sigma_{k+1}"] & S_{k+1}(\mathcal O^3) \arrow[d, ""] \arrow[r, "\tau_k"] & \mathfrak{K}_k\ar[r]&0\\
  0 \arrow[r] & R_{k+1} \arrow[r, ""] \ar[d, "\overline{\pi}_{k}"]& J_{k+3}(\mathcal{O}) \arrow[r, "P_{k+1}"] \arrow[d, "\pi_{k+2}"]& J_{k+1}(\mathcal{O}^3)\ar[d, "\hat{\pi}_k"] \\
  0 \arrow[r] & R_{k}\arrow[r, ""] \ar[d, "\nu_k"]& J_{k+2}(\mathcal O) \arrow[r, "P_k"] \arrow[d]& J_{k}(\mathcal{O}^3) \ar[d] \\
  &   \mathfrak{K}_k & 0&0
\end{tikzcd}
\end{equation}
$\mathfrak{K}_k:=\coker(\sigma_{k+1})$ is called the obstruction to formal integrability. Furthermore it is clear that $\overline{\pi}_{k}$ is surjective if and only if $\nu_k=0$. For later use we recall that for $m\geqslant1$ and $k\geqslant0$ we have $\rk J_k(\mathcal{O}^m)=\dfrac{m}{2}(k+1)(k+2)$ and $\rk S_k(\mathcal{O}^m)=m(k+1)$. We set $J_l(\mathcal{O}^m)=0$ for $l<0$. When $k\geqslant1$ $R_{k}$ and $g_k$ are $\mathcal{O}$-modules of finite type but not necessarily free. Also we have a natural exact sequence \cite[eq.~1.3.3]{spencer1969}
\begin{equation}
\label{abel404}
\begin{tikzcd}
  0 \arrow[r]& \mathcal{O}^m \arrow[r,"j_l"]& J_l(\mathcal{O}^m)\arrow[r, "D" ] & \Omega^1\otimes_\mathcal{O}J_{l-1}(\mathcal{O}^m)  \arrow[r, "D" ] &  \Omega^2\otimes_{\mathcal O}J_{l-2}(\mathcal{O}^m)\rar&0 \end{tikzcd}\end{equation}
 which stops at $\Omega^2\otimes_{\mathcal O}J_{l-2}(\mathcal{O}^m)$ due to the fact that our underlying space has dimension $2$. Let us compute the action of the map $D$ in local coordinates when $m=1$ and $l=2$. With previous notations if $\pi_0:J_1(\mathcal O)\longrightarrow J_0(\mathcal{O})$ is the natural map (projection) we have explicitly the formulas
 \begin{equation}
 \label{abel405}
 \begin{split}
 &D(z,p,q,r,s,t)\\
 &=dx\otimes(\partial_x(z)-p,\partial_x(p)-r,\partial_x(q)-s)\\
 &+du\otimes(\partial_u(z)-q,\partial_u(p)-s,\partial_u(q)-t)\\
 &D(\omega\otimes(z,p,q))\\
 &=d\omega\otimes\pi_0(z,p,q)-\omega\wedge D(z,p,q)\\
 &=d\omega\otimes (z)-\omega\wedge(dx\otimes(\partial_x(z)-p)+du\otimes(\partial_u(z)-q)).
 \end{split}
 \end{equation}
 Let us denote by $\mathscr{N}$ the system \eqref{abel40}. With notations as above the exact sequence \eqref{abel404} induces a complex on kernels called the first Spencer complex associated with the prolongation $p_k$ of $p_0$, that is the complex of (families) of $\C$-vector spaces exact at $R_k$ and with injective $j_{k+2}$ \cite[eq.~1.6.2, p.~196]{spencer1969}:
 \begin{equation}
 \label{abel406}
 \begin{tikzcd}
  0 \arrow[r]& \sol(\mathscr{N}) \arrow[r,"j_{k+2}"]& R_k\arrow[r, "D" ] & \Omega^1\otimes_\mathcal{O}R_{k-1}  \arrow[r, "D" ] &  \Omega^2\otimes_{\mathcal O}R_{k-2}\rar&0 \end{tikzcd}.
  \end{equation}
In equation \eqref{abel406} $\sol(\mathscr{N})$ is the set of $z\in\mathcal{O}$ whose jet of order $k+2$ lies in $R_k$ (another formulation of the notion of solution in the classical sense). Hence  we obtain the following commutative diagram of $\C$-vector spaces, from previous constructions. Its rows are exact and its columns are complexes exact at $R_k$ (resp. $R_{k+1}$) with injective $j_{k+2}$ (resp. $j_{k+3}$)
  \begin{equation}
  \label{abel407}
  \begin{tikzcd}
  &   & 0\ar[d] &0\ar[d]\\
 & 0 \arrow[r] & \sol(\mathscr N) \arrow[d, "j_{k+3}"] \arrow[r, "\Id"] & \sol(\mathscr N) \arrow[d, "j_{k+2}"] \arrow[r, ""] &  0\\
&  g_{k+1} \arrow[r] & R_{k+1} \arrow[r, "\overline{\pi}_{k}"] \ar[d, "D"]& R_{k} \arrow[r, "\nu_k"] \arrow[d, "D"]& \mathfrak{K}_k\\
 0\arrow[r]& \Omega^1\otimes_{\mathcal O}g_k \arrow[r] & \Omega^1\otimes_{\mathcal O}R_{k}\arrow[r, "\Id\otimes\overline{\pi}_{k-1}"] \ar[d, "D"]& \Omega^1\otimes_{\mathcal O}R_{k-1} \arrow[r, "\Id\otimes\nu_{k-1}"] \arrow[d,"D"]& \Omega^1\otimes_{\mathcal O}\mathfrak{K}_{k-1} \\
 0\arrow[r]& \Omega^2\otimes_{\mathcal O}g_{k-1} \arrow[r,""] & \Omega^2\otimes_{\mathcal O}R_{k-1}\arrow[r, "\Id\otimes\overline{\pi}_{k-2}"] \ar[d, ""]& \Omega^2\otimes_{\mathcal O}R_{k-2}\arrow[r, "\Id\otimes\nu_{k-2}"] \arrow[d]& \Omega^2\otimes_{\mathcal O}\mathfrak{K}_{k-2} \\
 & &   0 & 0&
  \end{tikzcd}
  \end{equation}
  where the last two rows of the diagram are due to the fact that $\Omega^1$ and $\Omega^2$ are free over $\mathcal{O}$ and to the first column of equation \eqref{abel403}.

The operator $D:R_k\to\Omega^1\otimes_\mathcal{O}R_{k-1}$ of the first Spencer sequence \eqref{abel404} satisfies
  $$D(fe)=df\otimes\overline{\pi}_{k-1}(e)+fD(e),\,f\in\mathcal{O},e\in R_k.$$
  Hence if $R_k$ is free and $\overline{\pi}_{k-1}$ is an isomorphism, i.e. $\nu_{k-1}=0$ and $g_k=0$ (equation \eqref{abel403}) this allows us to endow $R_k$ with a connection defined by
  $$\nabla_k:=\Id_{\Omega^1}\otimes\overline{\pi}_{k-1}^{-1}\circ D:R_k\to\Omega^1\otimes_{\mathcal O}R_k.$$
 Moreover we see that $\ker\nabla_k\simeq\sol(\mathscr N)$ since $\Id_{\Omega^1}\otimes\overline{\pi}_{k-1}^{-1}$ is an isomorphism (by flatness of $\Omega^1$ over $\mathcal O$) and $\ker D=\sol(\mathscr N)$; by diagram chasing in \eqref{abel407} we find that the curvature of $\nabla_k$, that is $(\Id_{\Omega^2}\otimes_{\mathcal O}\overline{\pi}_{k-1}^{-1}\circ D)\circ(\Id_{\Omega^1}\otimes_{\mathcal O}\overline{\pi}_{k-1}^{-1}\circ D)$ takes values in $\Omega^2\otimes_{\mathcal O}g_{k-1}$, up to isomorphism. Indeed we have
 $$(\Id_{\Omega^2}\otimes\overline{\pi}_{k-2})\circ D\circ(\Id_{\Omega^1}\otimes\overline{\pi}_{k-1})^{-1}\circ D=D\circ D=0.$$
   Hence $\image (D\circ(\Id_{\Omega^1}\otimes\overline{\pi}_{k-1})^{-1}\circ D)\subset \ker(\Id_{\Omega^2}\otimes \overline{\pi}_{k-2})=\Omega^2\otimes_{\mathcal O}g_{k-1} $. Therefore
$$(\Id_{\Omega^2}\otimes\overline{\pi}_{k-1})^{-1}D\circ(\Id_{\Omega^1}\otimes\overline{\pi}_{k-1})^{-1}\circ D$$   
takes values in $\Id_{\Omega^2}\otimes\overline{\pi}_{k-1}^{-1}(\Omega^2\otimes_{\mathcal O}g_{k-1})\simeq\Omega^2\otimes_{\mathcal O}g_{k-1}$ by flatness of $\Omega^2$ over $\mathcal O$.

Let us now compute the data for $P_0$. The kernel $R_0$ of $P_0$ is a free $\mathcal{O}$-module of rank 3 generated by
\begin{equation}
\label{abel401}
\begin{split}
&e_1=\left(\begin{array}{cccccc}1 & 0 & 0 & -r_3 & -r_2 & -r_1\end{array}\right)\\
&e_2=\left(\begin{array}{cccccc}0 & 1 & 0 & -p_3 & -p_2 & -p_1\end{array}\right)\\
&e_3=\left(\begin{array}{cccccc}0 & 0 & 1 & -q_3 & -q_2 & -q_1\end{array}\right).
\end{split}
\end{equation}
Also
   \begin{equation*}
   \begin{split}
   \pi_1\colon&J_2(\mathcal O)\to J_1(\mathcal O)\\
   &(z,p,q,r,s,t)\mapsto(z,p,q)
   \end{split}
   \end{equation*}
   hence $S_2(\mathcal O)=\ker(\pi_1)=\{(r,s,t),\,r,s,t\in\mathcal{O}\}$. This gives the symbol map
   \begin{equation*}
   \begin{split}
   \sigma_0\colon&S_2(\mathcal O)\to\mathcal{O}^3\\
   &(r,s,t)\mapsto(r,s,t)
   \end{split}
   \end{equation*}
   thus $g_0:=\ker(\sigma_0)=0$. Furthermore since $\sigma_0$ is surjective we have $\coker(\sigma_0)=0$. This gives $\mathfrak{K}_{-1}=0$ and an isomorphism
   \[\overline{\pi}_{-1}: R_0 \xrightarrow{\sim} R_{-1},\]
   where $R_{-1}=J_1(\mathcal O)$, as follows from equation \eqref{abel403}.
   
   Finally $g_{-1}$ is of rank $2$ because $S_{-1}(\mathcal O^3)=0$, since $J_{-1}(\mathcal O^3)=0$, which gives $S_1(\mathcal O)\simeq g_{-1}$ from the exact sequence \eqref{abel402}, and also because as recalled above $\rk S_1(\mathcal O)=2$. In conclusion we have the existence of a non-necessarily flat connection $\nabla_0$ on $R_0$ with curvature taking values in $\Omega^2\otimes_{\mathcal O}g_{-1}$. Let us determine its connection and curvature matrix in the basis $(e_1,e_2,e_3).$ We have by \eqref{abel405}
   \begin{equation}
   \label{abel4010000}
   \begin{split}
   &D(e_1)=dx\otimes(0,r_3,r_2)+du\otimes(0,r_2,r_1)\\
   &D(e_2)=dx\otimes(-1,p_3,p_2)+du\otimes(0,p_2,p_1)\\
   &D(e_3)=dx\otimes(0,q_3,q_2)+du\otimes(-1,q_2,q_1).
   \end{split}
   \end{equation}
   Computing $(\nabla_0(e_i))_{1\leqslant i\leqslant3}=(\Id_{\Omega^1}\otimes\overline{\pi}_{-1}^{-1}D(e_i))_{1\leqslant i\leqslant3}$ we obtain the following connection matrix in the adapted basis $(e_1,e_2,e_3)$
   \begin{equation}
   \label{abel40010001}
   \gamma=\left(\begin{array}{ccc}0 & -dx & -du \\r_3dx+r_2du & p_3dx+p_2du & q_3dx+q_2du \\r_2dx+r_1du & p_2dx+p_1du & q_2dx+q_1du\end{array}\right).\end{equation}
   Its curvature matrix is then given by
   \begin{equation}
   \label{abel40010003}
   d\gamma+\gamma\wedge\gamma=\left(\begin{array}{ccc}0 & 0 & 0 \\A_{21}& A_{22} & A_{23} \\A_{31} & A_{32} & A_{33}\end{array}\right)dx\wedge du\end{equation}
   with
   \begin{equation}
   \label{abel40010004}
   \begin{split}
   &A_{21}=-\partial_u(r_3)+\partial_x(r_2)+p_3r_2-p_2r_3+q_3r_1-q_2r_2 ,\\ 
   &A_{22}=-\partial_u(p_3)+\partial_x(p_2)+r_2+q_3p_1-q_2p_2\\
   &A_{23}=-\partial_u(q_3)+\partial_x(q_2)-r_3+p_3q_2-p_2q_3+q_3q_1-q_2^2\\
   &A_{31}=-\partial_u(r_2)+\partial_x(r_1)+p_2r_2-p_1r_3+q_2r_1-q_1r_2\\
   &A_{32}=-\partial_u(p_2)+\partial_x(p_1)+r_1+p_2^2-p_1p_3+q_2p_1-p_2q_1\\
   &A_{33}=-\partial_u(q_2)+\partial_x(q_1)-r_2+p_2q_2-p_1q_3.
   \end{split}
   \end{equation}
   Thus one sees that the connection $\nabla_0$ is flat, that is $\nabla_0\circ\nabla_0=0$, or again $d\gamma+\gamma\wedge\gamma=0$ if and only if the system of partial differential equations \eqref{abel40} satisfies all its six integrability conditions. Hence in this case the sheaf of horizontal sections of the connection $\nabla_0$: $\ker(\nabla_0)$ is a local system of rank three, \cite{deligne1970}.
\end{proof}

 Let us consider a second order differential equation of the form
\begin{equation}
\label{abel5}
u^{\prime\prime}+a_1(x,u)u^{\prime3}+3a_2(x,u)u^{\prime2}+3a_3(x,u)u^{\prime}+a_4(x,u)=0
\end{equation}
with $a_1(x,u)$, $a_2(x,u)$, $a_3(x,u)$, $a_4(x,u)\in \C\{x,u\}$, which defines a projective connection. Then the necessary and sufficient condition for the existence of a point transformation sending equation \eqref{abel5} to $u^{\prime\prime}=0$ is \cite{cartan1924, liouville18892, tresse1894, olver}
\begin{equation}
\label{abel6}
\begin{split}
L_1&=\dfrac{\partial}{\partial x}\left(2\dfrac{\partial a_3}{\partial u}-\dfrac{\partial a_2}{\partial x}+a_1a_4\right)-\dfrac{\partial}{\partial u}\left(\dfrac{\partial a_4}{\partial u}+3a_2a_4\right)+a_4\left(\dfrac{\partial a_1}{\partial x}+3a_1a_3\right)\\
&+3a_3\left(2\dfrac{\partial a_3}{\partial u}-\dfrac{\partial a_2}{\partial x}+a_1a_4\right)=0\\
L_2&=\dfrac{\partial}{\partial u}\left(\dfrac{\partial a_3}{\partial u}-2\dfrac{\partial a_2}{\partial x}+a_1a_4\right)+\dfrac{\partial}{\partial x}\left(\dfrac{\partial a_1}{\partial x}-3a_1a_3\right)+a_1\left(\dfrac{\partial a_4}{\partial u}+3a_2a_4\right)\\
&-3a_2\left(\dfrac{\partial a_3}{\partial u}-2\dfrac{\partial a_2}{\partial x}+a_1a_4\right)=0.
\end{split}
\end{equation}
The conditions \eqref{abel6} are equivalent to the integrability conditions of the following system of partial differential equations
\begin{equation*}
\begin{split}
z^{(2,0)}+a_3z^{(1,0)}-a_4z^{(0,1)}&+(\frac{\partial a_4}{\partial u}-\frac{\partial a_3}{\partial x}+2(a_2a_4-a_3^2))z=0\\
z^{(1,1)}+a_2z^{(1,0)}-a_3z^{(0,1)}&+(\frac{\partial a_3}{\partial u}-\frac{\partial a_2}{\partial x}+(a_1a_4-a_2a_3))z=0\\
z^{(0,2)}+a_1z^{(1,0)}-a_2z^{(0,1)}&+(\frac{\partial a_2}{\partial u}-\frac{\partial a_1}{\partial x}+2(a_1a_3-a_2^2))z=0\\
z^{(i,j)}&=\dfrac{\partial^{i+j}z}{\partial x^i\partial u^j},\,0\leqslant i,j\leqslant2
\end{split}
\end{equation*}
that is $z^{(1,2)}=\dfrac{\partial}{\partial x}(z^{(0,2)})=\dfrac{\partial}{\partial u}(z^{(1,1)})$ and $z^{(2,1)}=\dfrac{\partial}{\partial u}(z^{(2,0)})=\dfrac{\partial}{\partial x}(z^{(1,1)})$ are equivalent to
$$L_1=L_2=0$$
i.e. equation \eqref{abel6}, see \cite{liouville18872}, \cite[p.~532]{henaut1993}. Indeed in the notations of the \theoref{sec1} and its proof (see equations \eqref{abel40}, \eqref{abel40010003}, \eqref{abel40010004}) let us show that $A_{22}$, $A_{23}$, $A_{32}$, $A_{33}$ vanish. Then $A_{21}=0$, resp. $A_{31}=0$ will give after differentiation and elimination $L_1=0$, resp. $L_2=0$. We have $p_3=a_3$, $q_3=-a_4$, $p_2=a_2$, $q_2=-a_3$, $p_1=a_1$, $q_1=-a_2$, and
\begin{equation*}
\begin{split}
r_3&=\frac{\partial a_4}{\partial u}-\frac{\partial a_3}{\partial x}+2(a_2a_4-a_3^2)\\
r_2&=\frac{\partial a_3}{\partial u}-\frac{\partial a_2}{\partial x}+(a_1a_4-a_2a_3)\\
r_1&=\frac{\partial a_2}{\partial u}-\frac{\partial a_1}{\partial x}+2(a_1a_3-a_2^2).
\end{split}
\end{equation*}
This gives
\begin{equation}
\begin{split}
A_{22}&=-\partial_u(a_3)+\partial_x(a_2)+(-a_1a_4+a_2a_2)+r_2=0\\
A_{23}&=\partial_u(a_4)-\partial_x(a_3)-a_3^2+a_2a_4+a_2a_4-a_3^2-r_3=0\\
A_{32}&=-\partial_u(a_2)+\partial_x(a_1)+a_2^2-a_1a_3-a_1a_3+a_2^2+r_1=0\\
A_{33}&=\partial_u(a_3)-\partial_x(a_2)-a_2a_3+a_1a_4-r_2=0.
\end{split}
\end{equation}
\begin{corollary}[\cite{liouville18892}]
\label{lemma4}
The following system of linear partial differential equations
 \begin{equation}
 \label{abel16}
\begin{split}
z^{(2,0)}+a_3z^{(1,0)}-a_4z^{(0,1)}&+(\frac{\partial a_4}{\partial u}-\frac{\partial a_3}{\partial x}+2(a_2a_4-a_3^2))z=0\\
z^{(1,1)}+a_2z^{(1,0)}-a_3z^{(0,1)}&+(\frac{\partial a_3}{\partial u}-\frac{\partial a_2}{\partial x}+(a_1a_4-a_2a_3))z=0\\
z^{(0,2)}+a_1z^{(1,0)}-a_2z^{(0,1)}&+(\frac{\partial a_2}{\partial u}-\frac{\partial a_1}{\partial x}+2(a_1a_3-a_2^2))z=0\\
L_1=L_2=0,\,z^{(i,j)}&=\dfrac{\partial^{i+j}z}{\partial x^i\partial u^j},\,0\leqslant i,j\leqslant2\end{split}
\end{equation}
has three linearly independent solutions in $\C\{x,u\}$ over $\C$: the dimension of its solution space $\mathscr{S}\subseteq\C\{x,u\}$ over $\C$ is three, where we recall that $L_1=L_2=0$ are equivalent to the integrability of the system \eqref{abel16}. 
\end{corollary}
\begin{proof}
This immediately follows from \theoref{sec1}. Indeed since the system \eqref{abel6}: $L_1=L_2=0$, is satisfied, the integrability conditions of the system \eqref{abel16} are all verified. Hence the curvature of the connection provided by \theoref{sec1} vanishes, thus the solution space of the system has dimension $3$ over $\C$.
 \end{proof}
Let us see how to get from the system \eqref{abel16} to the second order ordinary differential equation
\eqref{abel5}.
\begin{lemma}
\label{lemma5}
Let $z_1$, $z_2$, $z_3\in\C\{x,u\}$ be three linear independent solutions to equation \eqref{abel16} and $z=\alpha_1z_1+\alpha_2z_2+\alpha_3z_3$, $\alpha_1$, $\alpha_2$, $\alpha_3\in\C$ its general solution then by setting $z=0$ (which defines implicitly $u$ in terms of $x$) and using the system \eqref{abel16} we obtain the differential equation \eqref{abel5}
$$u^{\prime\prime}+a_1(x,u)u^{\prime3}+3a_2(x,u)u^{\prime2}+3a_3(x,u)u^{\prime}+a_4(x,u)=0.$$
\end{lemma}
\begin{proof}
Since $z(x,u)=\alpha_1z_1+\alpha_2z_2+\alpha_3z_3$ is the general solution to equation \eqref{abel16} and $z_1$, $z_2$, $z_3$ are linearly independent, we have $z^{(0,1)}=z_x=\dfrac{\partial}{\partial x}(z)\not=0$. Thus we can apply the implicit function theorem to define a function $u(x)$ from the equation $z(x,u)$=0. As
$$z(x,u(x))=0$$
we deduce
 \begin{equation}
 \label{abel17}
\begin{split}
z^{(0,1)}u^\prime+z^{(1,0)}=0\\
z^{(2,0)}+2z^{(1,1)}u^\prime+z^{(0,2)}u^{\prime2}&+z^{(0,1)}u^{\prime\prime}=0.
\end{split}
\end{equation}
We next use the system \eqref{abel16} and the first equation of \eqref{abel17} to replace the expressions $z^{(1,0)}$, $z^{(1,1)}$, $z^{(2,0)}$, $z^{(0,2)}$ in the second equation \eqref{abel17} by their values in terms of $z^{(0,1)}$, $z$ $u^\prime$ and the $a_i(x,u)_{1\leqslant i\leqslant4}$ and their derivatives. We get a relation of the following form
 \begin{equation*}
 \label{abel18}
\begin{split}
&(u^{\prime\prime}+a_1(x,u)u^{\prime3}+3a_2(x,u)u^{\prime2}+3a_3(x,u)u^{\prime}+a_4(x,u))z^{(0,1)}(x,u(x))+\\&
g(x,u(x))z(x,u(x))=0, \;g(x,u)\in\mathcal{O}.
\end{split}
\end{equation*}
This gives
 \begin{equation}
 \label{abel1800}
(u^{\prime\prime}+a_1(x,u)u^{\prime3}+3a_2(x,u)u^{\prime2}+3a_3(x,u)u^{\prime}+a_4(x,u))z^{(0,1)}(x,u(x))=0.
\end{equation}
Since $z^{(0,1)}(x,u(x))\not=0$ we see that equation \eqref{abel1800} is equivalent to
$$u^{\prime\prime}+a_1(x,u)u^{\prime3}+3a_2(x,u)u^{\prime2}+3a_3(x,u)u^{\prime}+a_4(x,u)=0.$$
\end{proof}
\section{Characterization of some projective connections}
\label{scheffers}
In this section we study more in depth some second order differential equations of projective connections that is differential equations of the form
\begin{equation}
\label{pro11}
\begin{split}
&u^{\prime\prime}=A(x,u)+B(x,u)u^\prime+C(x,u)u^{\prime2}+D(x,u)u^{\prime3},\\ 
&A(x,u), B(x,u), C(x,u), D(x,u)\in \mathcal{Q}:=\FRAC(\mathcal O).
\end{split}
\end{equation}
Let $P=(x_P,u_P)$ belong to the plane $(\C^2,0)$. We consider on the plane $(\C^2,0)$ three germs of lines $\mathscr L_1$, $\mathscr L_2$, $\mathscr L_3$ with respective directions $w_1$, $w_2$, $w_3$ and which pass through the point $P$. Further let us take three points $Q_1\in \mathscr L_1$, $Q_2\in \mathscr L_2$ and $Q_3\in \mathscr L_3$ (below in this section, we are interested in the "generic" case where no three of the points $P$, $Q_1$, $Q_2$, $Q_3$, lie on a line). Then the four points $P$, $Q_1$, $Q_2$, $Q_3$ define a  germ of pencil of conics which contain them and every member of this pencil is uniquely determined once we fix the direction $u^\prime=\frac{du}{dx}$ of its tangent at $P$ (reasoning similar to the case where one has to show that there is a unique conic through five given points, we refer to \cite[chap.~7]{semplekneebone}, for more on pencil of conics). Furthermore if the tangent line at $P$ is required to pass through any of $Q_1$, $Q_2$, $Q_3$, then the conic will be degenerate or singular. By fixing the direction of the tangent at $P$, one also determines the value of $u^{\prime\prime}=\frac{d^2u}{dx^2}$ at $P$, since once the conic is uniquely fixed, the value of $u^{\prime\prime}$ at $P$ is also uniquely fixed. Moreover this implies that $u^{\prime\prime}$ is a polynomial function of $u^\prime$. Indeed the general equation of a conic through $P$ is
\begin{equation}
\label{oooo}
\begin{split}
a_0x^2+2b_0xu+c_0u^2+2d_0x+2e_0u+f_0=0.\end{split}
\end{equation}
If $c_0\not=0$ then $u$ can be expressed explicitly in the form
$$u=m_0x+n_0+(p_0x^2+2q_0x+r_0)^{1/2},$$
for $m_0$, $n_0$, $p_0$, $q_0$, $r_0\in\C$. If $p_0=q_0=r_0$ then $u^{\prime\prime}=0$. Otherwise we obtain
$$u^{\prime\prime}=(p_0r_0-q_0^2)(p_0x^2+2q_0x+r_0)^{-3/2}.$$
Besides
$$u^{\prime}=m_0+\frac{p_0x+q_0}{(p_0x^2+2q_0x+r_0)^{1/2}}.$$
If $p_0=q_0=0$ then $u^{\prime\prime}=0$. Otherwise we have
$$\frac{1}{(p_0x^2+2q_0x+r_0)^{1/2}}=\frac{u^\prime-m_0}{p_0x+q_0}.$$
Hence $u^{\prime\prime}$ is a polynomial in $u^\prime$ of degree $3$ with $\mathcal Q$ coefficients, in a sufficiently small neighborhood of $P$. If $c_0=0$ then $u$ has the following explicit form
$$u=\alpha x+\beta+\dfrac{\gamma}{b_0x+e_0},\;\alpha,\beta,\gamma\in\C,$$
and $b_0$, $e_0$ can not both be zero, by equation \eqref{oooo}. In case $\gamma=0$ then $u^{\prime\prime}=0$. If $\gamma\not=0$ then
$$u^{\prime\prime}=\dfrac{2\gamma b_0^2}{(b_0x+e_0)^3}=:A_0(x,u),$$
where $A_0(x,u)\in\mathcal Q$.

Let us designate by $(Q_1Q_2)$, resp. $(Q_1Q_3)$, resp. $(Q_2Q_3)$, the line through $Q_1$ and $Q_2$, resp. $Q_1$ and $Q_3$, resp. $Q_2$ and $Q_3$. By construction we can write the equation of the pencil of conics through the four points $P$, $Q_1$, $Q_2$, $Q_3$, in either of the following forms 
\begin{equation*}
\begin{split}
&\alpha_1\mathscr L_1(Q_2Q_3)+\alpha_2\mathscr L_2(Q_1Q_3)=0,\quad\beta_1\mathscr L_1(Q_2Q_3)+\beta_2\mathscr L_3(Q_1Q_2)=0,\\
&\gamma_1\mathscr L_2(Q_1Q_3)+\gamma_2\mathscr L_3(Q_1Q_2)=0 
\end{split}
\end{equation*}
where $\left[\alpha_1:\alpha_2\right]$, $\left[\beta_1:\beta_2\right]$, $\left[\gamma_1:\gamma_2\right]$, are arbitrary points in $\C P_1$.
 Now the pencil of conics through the four given points is degenerate if and only if, exactly one of $\alpha_1$, $\alpha_2$, or $\beta_1$, $\beta_2$, or $\gamma_1$, $\gamma_2$, is zero. Thus if the pencil is degenerate, $u^\prime$ takes at $P$ either of the values $w_1$, $w_2$, $w_3$. Conversely if $u^\prime$ takes at $P$ either of the values $w_1$, $w_2$, $w_3$, then by unicity, the conic is either of the following product of lines
 \begin{equation*}
 \begin{split}
 &\mathscr L_1(Q_2Q_3)=0, \; \mathscr L_2(Q_1Q_3)=0, \;\mathscr L_1(Q_2Q_3)=0, \\
 &\mathscr L_3(Q_1Q_2)=0,\; \mathscr L_2(Q_1Q_3),\;  \mathscr L_3(Q_1Q_2)=0,
 \end{split}
 \end{equation*}
 that is the conic is degenerate. Thus degeneracy of the pencil is equivalent to $u^{\prime\prime}=0$ at $P$. Therefore we find, in the "generic" case that 
\begin{equation}
\label{pro12}
u^{\prime\prime}=D_0(u^\prime-w_1)(u^\prime-w_2)(u^\prime-w_3)
\end{equation}
where $D_0$ does not involve $u^\prime$. Now if we vary $Q_1$, $Q_2$, $Q_3$ and allow $u^\prime$ at $P$ to remain the same, then the pencil of conics that we obtain will verify at $P$ a relation similar to \eqref{pro12}, given by
\begin{equation}
u^{\prime\prime}=D_1(u^\prime-w_1)(u^\prime-w_2)(u^\prime-w_3)
\end{equation}
but with $D_0$ necessarily different from $D_1$, because otherwise the two pencils would be identical: they would have the same value of $u^\prime$ at $P$ and the same value of $u^{\prime\prime}$ at $P$. Since we can furthermore vary the directions $w_1$, $w_2$, $w_3$, we obtain through this construction, the family of all equations of the form
\begin{equation*}
\begin{split}
u^{\prime\prime}&=D(x,u)(u^\prime-w_1(x,u))(u^\prime-w_2(x,u))(u^\prime-w_3(x,u))\\
&=:A(x,u)+B(x,u)u^\prime+C(x,u)u^{\prime2}+D(x,u)u^{\prime3}\end{split}
\end{equation*}
$A$, $B$, $C$, $D\in\mathcal Q$, $w_1$, $w_2$, $w_3\in\mathcal{O}$. Thus
\begin{lemma}
\label{lem1}
Given an equation
\begin{equation}
\label{pro111}
\begin{split}
u^{\prime\prime}&=D(x,u)(u^\prime-w_1(x,u))(u^\prime-w_2(x,u))(u^\prime-w_3(x,u))\\
&=:A(x,u)+B(x,u)u^\prime+C(x,u)u^{\prime2}+D(x,u)u^{\prime3},\;A, B, C, D\in\mathcal Q,\;w_1,w_2,w_3\in\mathcal{O},
\end{split}
\end{equation}
and $P=(x_P,u_P)\in(\C^2,0)$, a non-singular point of the coefficients $A$, $B$, $C$, $D$, then there exists a germ of pencil of conics in $(\C^2,0)$, which contains $P$ and such that the corresponding values of the derivatives $u^\prime$, $u^{\prime\prime}$ with respect to $x$ at the point $P$ satisfy the given equation
$$u^{\prime\prime}=A(x,u)+B(x,u)u^\prime+C(x,u)u^{\prime2}+D(x,u)u^{\prime3}.$$
\end{lemma}
\definecolor{ududff}{rgb}{0.30196078431372547,0.30196078431372547,1}
\definecolor{xdxdff}{rgb}{0.49019607843137253,0.49019607843137253,1}
\definecolor{cqcqcq}{rgb}{0.7529411764705882,0.7529411764705882,0.7529411764705882}
\begin{center}
\begin{tikzpicture}[line cap=round,line join=round,>=triangle 45,x=1cm,y=1cm]
\draw [rotate around={-178.57164469985398:(-8.987038372736885,0.09754670059542438)},line width=0.5pt] (-8.987038372736885,0.09754670059542438) ellipse (3.0149859953691514cm and 1.196459769388341cm);
\draw [rotate around={13.776720048286453:(-9.183147026906486,-0.13182312385578993)},line width=0.5pt] (-9.183147026906486,-0.13182312385578993) ellipse (3.710225334363145cm and 1.1668458143237868cm);
\draw [rotate around={-14.68420828004502:(-9.795693117921148,0.5309953255755083)},line width=0.5pt] (-9.795693117921148,0.5309953255755083) ellipse (4.0882156575533735cm and 1.229168069726762cm);
\begin{scriptsize}
\draw [fill=xdxdff] (-12,0) circle (2.5pt);
\draw[color=ududff] (-12,0.458456483989531) node {$P$};
\draw [fill=ududff] (-6,0) circle (2.5pt);
\draw[color=ududff] (-5.82495382672716,0.458456483989531) node {$Q_1$};
\draw [fill=ududff] (-8.049124944093071,1.2546086453193734) circle (2.5pt);
\draw[color=ududff] (-7.872202241575328,1.7221900734019793) node {$Q_2$};
\draw [fill=ududff] (-8.049124944093071,-1.0201118156230338) circle (2.5pt);
\draw[color=ududff] (-7.872202241575328,-0.5525303875404277) node {$Q_3$};
 \draw [red](-15,0) -- (-2,0);
  \draw (-12,0) -- (-8.05,1.25);
  \draw[red](-15,-0.95) -- (-2,3.16);
  \draw [red](-15,0.77) -- (-2,-2.532);
      \node [black] at (-1.7, 0) {$\mathscr L_1$};
      \node [black] at (-2, 3.4) {$\mathscr L_2$};
      \node [black] at (-1.7, -2.6) {$\mathscr L_3$};

\end{scriptsize}
\end{tikzpicture}
 \captionof{figure}{:\,Visualization of some of the members of the pencil of conics through $P$, $Q_1$, $Q_2$, $Q_3$, together with the three lines $\mathscr L_1$, $\mathscr L_2$, $\mathscr L_3$.}
\label{default}
\end{center}


\begin{definition}Consider a differential of the second order
$$u^{\prime\prime}=f(x,u,u^\prime)$$
with $f\in\mathcal Q\left[u^\prime\right]$. We say that the set of its integral curves in the plane $(\C^2,0)$ possesses a pencil of second order contact conics at a non-singular point $P\in(\C^2,0)$ of its coefficients, if and only if there exists a germ of pencil of conics which passes through $P$, and such that every germ of integral curve of $u^{\prime\prime}=f(u,u,u^\prime)$ which passes through $P$ has a contact of second order with a conic of the pencil. This means that if we parameterize the germ of integral curve of $u^{\prime\prime}=f(u,u,u^\prime)$ by $(x,u(x), u^\prime(x))$, $u(x)\in\C\{x\}$, there exists a conic of the pencil: $\mathfrak g(x,u)=0$ such that
$$\mathfrak g(x,u(x))=(\mathfrak g(x,u(x)))^\prime=(\mathfrak g(x,u(x)))^{\prime\prime}=0,$$
where $^\prime=\frac{d}{dx}$.
\end{definition}
We have the
\begin{theorem}
\label{theo3}
The differential equations of the form \eqref{pro111}
are the unique second order differential equations $u^{\prime\prime}=f(x,u,u^\prime)$, $f\in\mathcal Q\left[u^\prime\right]$, such that the set of their integral curves in the plane $(\C^2,0)$, possesses at each non-singular point of its coefficients (which belong to $\mathcal Q$), a pencil of second order contact conics.
\end{theorem}
\begin{proof}
Suppose that a germ of integral curve of a second order differential equation $u^{\prime\prime}=f(x,u,u^\prime)$, or again a trajectory of the vector field $\frac{\partial}{\partial x}+u^\prime\frac{\partial}{\partial u}+f(x,u,u^\prime)\frac{\partial}{\partial u^\prime}$, has a second order contact at $P$ with a plane conic with constant complex coefficients of equation 
$$a_1x^2+b_1xu+c_1u^2+d_1x+e_1u+f_1=0.$$
This means that if we parameterize the $1$-jet of the integral curve by $(x,u(x),u^\prime(x))$ with $u(x)\in\C\{x\}$ then we have identically 
$$a_1x^2+b_1xu(x)+c_1u^2(x)+d_1x+e_1u(x)+f_1=0$$
as well as the following two equations gotten from the previous equation by derivation
\begin{equation}
\label{proj3009}
2a_1x+b_1u(x)+b_1xu^\prime(x)+2c_1uu^\prime+d_1+e_1u^\prime=0
\end{equation}
and
\begin{equation}
\label{proj4009}
2a_1+2b_1u^\prime+b_1xu^{\prime\prime}+2c_1u^{\prime2}+2c_1uu^{\prime\prime}+e_1u^{\prime\prime}=0.\end{equation}
If $e_1+2c_1u+b_1x\equiv0$ then necessarily $c_1\not=0$ since otherwise $e_1=b_1=0$ and then $a_1=f_1=d_1=0$, a contradiction. Hence if $e_1+2c_1u+b_1x\equiv0$, then $u^{\prime\prime}(x)=0$. Now if $e_1+2c_1u+b_1x\not\equiv0$, then from equation \eqref{proj3009} we deduce
$$e_1+2c_1u+b_1x=\frac{2a_1x+b_1u+d_1}{u^\prime}$$
where $u^\prime\not\equiv0$ (in case $u^\prime\equiv 0$ then $u^{\prime\prime}\equiv 0$).

Now inserting this found value of $e_1+2c_1u+b_1x$ into \eqref{proj4009}, we find after collecting that $u^{\prime\prime}=f(x,u,u^\prime)$ is a polynomial in $u^\prime$ of degree $3$ with coefficients in $\mathcal Q$, which is of the form \eqref{pro111}. Conversely according to \lemref{lem1} given a differential equation of a projective connection
$$u^{\prime\prime}=A(x,u)+B(x,u)u^\prime+C(x,u)u^{\prime2}+D(x,u)u^{\prime3},\,A, B, C, D\in\mathcal Q$$
of the form \eqref{pro111} and $P=(x_P,u_P)\in(\C^2,0)$, a non-singular point of its coefficients, one can find a germ of pencil of conics in $(\C^2,0)$ through $P$ such that the corresponding values of the derivatives $u^\prime$, $u^{\prime\prime}$ with respect to $x$ at the point $P$ verify the given differential equation. This implies that the germ of pencil of conics through $P$ has contact of second order at $P$ with the unique integral curve of \eqref{pro111} prescribed by the datum of $(P=(x_P,u_P), u^\prime(P))$.
\end{proof}
 \section{Differential Geometry of second order differential equations}
 \label{kasner}
 Let 
 \begin{equation}
 \label{pro1}
 \{X=\phi(x,u), U=\psi(x,u),\,\phi(x,u), \psi(x,u)\in\mathcal{O}:=\C\{x,u\},\,J:=\phi_x\psi_u-\phi_u\psi_x\not=0\}
 \end{equation}
 be the pseudo-group $\mathscr{P}$ of point transformations of the plane $(\C^2,0)$. 
 
 Let $u^\prime$ and $u^{\prime\prime}$ be jet coordinates on $J_2(\mathcal O)$, the second order jet space. We study in this section the action of $\mathscr{P}$ on the tuple $(u^\prime,u^{\prime\prime})$ and deduce some consequences on differential equations of the second order polynomial in $u^{\prime\prime}$, $u^\prime$, and with holomorphic coefficients in $\mathcal{O}$. 
 
 The action of $\mathscr{P}$ on $(u^\prime,u^{\prime\prime})$ is given by
 \begin{equation}
 \label{pro3}
 \begin{split}
 &U^\prime=\dfrac{\psi_x+u^\prime\psi_u}{\phi_x+u^\prime\phi_u}\\
 &U^{\prime\prime}=\dfrac{\lambda+\mu u^\prime+\nu u^{\prime2}+\xi u^{\prime3}+Ju^{\prime\prime}}{(\phi_x+u^\prime\phi_u)^3}.
 \end{split}
 \end{equation}
 Here $\lambda$, $\mu$, $\nu$, $\xi$ belong to $\mathcal{O}$ and depend on derivatives of $\phi$ and $\psi$ up to order $2$. Furthermore we may consider $\lambda$, $\mu$, $\nu$, $\xi$ as well as $a:=\psi_u$, $b:=\psi_x$, $c:=\phi_u$ and $d:=\phi_x$ as numbers since we are only interested in the variations of the tuple $(u^\prime,u^{\prime\prime})$. In order to simplify our expressions and calculations we introduce the notation
 $$v:=u^\prime,\,w=u^{\prime\prime};\,V:=U^\prime,W:=U^{\prime\prime}.$$
 Thus the transformation \eqref{pro3} takes the form
 \begin{equation}
 \label{pro4}
 \begin{split}
 V=\dfrac{av+b}{cv+d},\,W=\dfrac{\lambda+\mu v+\nu v^2+\xi v^3+(ad-bc)w}{(cv+d)^3}.
 \end{split}
 \end{equation}
 Hence the pseudo-group $\mathscr{P}$ induces on the tuple $(v,w)$ the eight-parameter (local) Lie group of transformations \eqref{pro4} whose infinitesimal generators are given by
 \begin{equation}
 \label{pro5}
 \dfrac{\partial}{\partial v},\,\dfrac{\partial}{\partial w},\,v\dfrac{\partial}{\partial v},\,w\dfrac{\partial}{\partial w},\,v^2\dfrac{\partial}{\partial v}+3vw\dfrac{\partial}{\partial w},\,v\dfrac{\partial}{\partial w},v^2\dfrac{\partial}{\partial w},\,v^3\dfrac{\partial}{\partial w}.
 \end{equation} 
 Together these vector fields generate a Lie algebra isomorphic to $\mathfrak{g}\mathfrak{l}(2,\C)\rtimes S^3\C^2$, the Lie algebra of the Lie group $GL(2,\C)\rtimes S^3\C^2$, where $S^3\C^2\simeq\C^4$ designates the third symmetric product of $\C^2$ with itself, see the number $28$ of the list given in \cite[p.~341]{olver1}, see also \cite[p.~472]{olver}. Here one interprets $S^3\C^2$ as the dual of $S^3((\C^2)^*)$, which is the space of complex homogeneous polynomials of degree $3$ in two variables.
 \subsection{Link with the four dimensional projective space $\C P_4$}We remark that the numerators and denominators of $V$ and $W$ are polynomials in $v$ and $w$ of degree at most $3$. Moreover they are linear in $v$, $v^2$, $v^3$, $w$. Thus we introduce the following point of $\C P_4$ given in terms of homogeneous coordinates by
 \begin{equation}
 \label{pro6}
 \left[z_1:z_2:z_3:z_4:z_5\right]:=\left[1:v:v^2:v^3:w\right].
 \end{equation}  
 This gives an embedding of $(v,w)$ in $\C P_4$. In $\C P_4$ $(v,w)$ is characterized by the surface $\mathscr S_3$ given by the following homogenous equations
 \begin{equation}
 \label{pro7}
 z_1z_3-z_2^2=0,\quad z_2z_3-z_1z_4=0,\quad z^2_3-z_2z_4=0.
 \end{equation} 
 $\mathscr S_3$ is cone in $\C P_4$ with vertex given by the point $\left[0:0:0:0:1\right]$. The section of $\mathscr S_3$ in the plane $z_5=0$ is a twisted cubic. Hence the subscript $3$ in $\mathscr S_3$. The transformations \eqref{pro4} may now be written in the following form
 \begin{equation}
 \label{pro8}
 \begin{split}
& Z_1=d^3z_1+3cd^2z_2+3c^2dz_3+c^3z_4\\
 &Z_2=b^2dz_1+d(ad+2bc)z_2+c(bc+2ad)z_3+ac^2z_4\\
 &Z_3=bd^2z_1+b(bc+2ad)z_2+a(ad+2bc)z_3+a^2cz_4\\
 &Z_4=b^3z_1+3ab^2z_2+3a^2bz_3+a^3z_4\\
 &Z_5=\lambda z_1+\mu z_2+\nu z_3+\xi z_4+(ad-bc)z_5.
 \end{split}
 \end{equation}
where $\left[Z_1:Z_2:Z_3:Z_4:Z_5\right]=\left[1:V:V^2:V^3:W\right]$. The transformations \eqref{pro8} leave the relations \eqref{pro7} invariant. Thus to each transformation \eqref{pro4} of $(u^\prime,u^{\prime\prime})=(v,w)$ there corresponds in $\C P_4$ an element of $PGL(4,\C)$ which leaves the cubic cone $\mathscr S_3$ invariant. We denote the set of all projective transformations of $\C P_4$ which arise in this way by $G$.

One can show that the dimension of a Lie group of projective transformations which leaves a cubic cone in $\C P_4$ invariant is at most $8$. Indeed the dimension of $PGL(4,\C)$, the group of automorphisms of $\C P_4$, is $24$. Let $H$ be the subgroup of $PGL(4,\C)$ in question. It has to fix the vertex $\left[0:0:0:0:1\right]$ of the cone. This reduces its dimension by $4$ and leaves us with $20$ as the possible dimension of $H$. We recall \cite{harris1992, harris1998}, that the dimension of the space of twisted cubics in $\C P_4$ is $16$. Indeed a twisted cubic in $\C P_4$ is the image of a holomorphic (polynomial) map
\begin{equation*}
\begin{split}
&\C P_1\to \C P_4\\
&\left[x_0:x_1\right]\mapsto \left[P_0(x_0,x_1),P_1(x_0,x_1),P_2(x_0,x_1),P_3(x_0,x_1), P_4(x_0,x_1)\right]
\end{split}
\end{equation*}
where $P_0(x_0,x_1),\ldots,P_4(x_0,x_1)$, are homogeneous polynomials of degree $3$ in $x_0$, $x_1$. Moreover these choices determine a twisted cubic up to scalar multiplication and automorphisms of $\C P_1$, that is $PGL(2,\C)$. Now $P_0(x_0,x_1),\ldots,P_4(x_0,x_1)$ bring $20$ parameters from which we have to subtract the dimension of $PGL(2,\C)$, that is $3$, and $1$ due to scalar multiplication. This gives the mentioned number $16$. Also the dimension of the space of twisted cubics on the cone over 
$$ z_1z_3-z_2^2=0,\quad z_2z_3-z_1z_4=0,\quad z^2_3-z_2z_4=0$$
is $4$. In fact such a twisted cubic on the cone is a point of the $16$-dimensional space of twisted on $\C P_4$ which belongs to the cone. Thus it is parameterized as follows
$$\left[x_0:x_1\right]\mapsto\left[x_0^3:x_0^2x_1:x_0x_1^2:x_1^3:P_5(x_0,x_1)\right]$$
with $P_5(x_0,x_1)$ a homogeneous polynomial in $x_0$, $x_1$ of degree 3. Therefore we get the stated dimension $4$.

Let us finish determining the dimension of $H$. It is required that the twisted cubic in the plane $z_5=0$ be converted into one of the $4$-dimensional space of twisted cubics on the cone. Since the dimension of the space of twisted cubics in $\C P_4$ is $16$ this imposes $16-4=12$ restrictions. Hence the group $H$ in question contains at most $20-12=8$ independent parameters hence precisely $8$ parameters. Thus $H$ is of dimension $8$ and $H=G$. Therefore
\begin{proposition}
\label{theo1}
The study of the geometry of $(u^\prime,u^{\prime\prime})$ under arbitrary point transformations of the plane is equivalent to the study of the projective geometry on a cubic cone in $\C P_4$.
\end{proposition} 
\subsection{Necessary condition for equivalence of certain second order ODEs}
\propref{theo1} suggests a necessary condition of equivalence for second order ordinary differential equations which are polynomial in $u^\prime$, $u^{\prime\prime}$, that is in $v$, $w$, according to the degree of the corresponding equation expressed in $z$ coordinates. To an equation of degree one 
$$A_1z_1+A_2z_2+A_3z_3+A_4z_4+A_5z_5=0$$
corresponds in terms of $v$, $w$ an equation
$$A_1+A_2v+A_3v^2+A_4v^3+A_4w=0.$$
Let us now take the coefficients as arbitrary functions of $x$, $u$ belonging to $\mathcal{O}$, then we have a differential equation of the form
\begin{equation}
\label{pro9}
Eu^{\prime\prime}=Du^{\prime3}+Cu^{\prime2}+Bu^{\prime}+A
\end{equation}
which defines a $2$-dimensional complex projective connection.

The next case, corresponding to equations of the second degree in the $z$ leads to an equation of the form
\begin{equation}
\label{pro10}
\begin{split}
&Au^{\prime\prime2}+(B_1+B_2u^\prime+B_3u^{\prime2}+B_4u^{\prime3})u^{\prime\prime}\\&+(C_0+C_1u^\prime+C_2u^{\prime2}+C_3u^{\prime3}+C_4u^{\prime4}+C_5u^{\prime5}+C_6u^{\prime6})=0
\end{split}
\end{equation}
where the coefficients belong to $\mathcal{O}$. In general we define the degree of a second order differential equation, polynomial in $u^\prime$, $u^{\prime\prime}$ and which can be written as the vanishing set of a homogeneous polynomial in the $(z_i)_{1\leqslant i\leqslant4}$ coordinates, as the degree in the $(z_i)_{1\leqslant i\leqslant4}$ of the corresponding equation in $(z_i)_{1\leqslant i\leqslant4}$ . 

The set of second order differential equations (for which one can define the degree) of any fixed degree $r$ is invariant under point transformations of the plane. So for two second order equations for which one can define their degrees, a necessary condition for their equivalence under point transformations of the plane, is the equality of their degrees. 

\begin{proposition}
\label{theo2}
Consider two second order differential equations polynomial in the first and second derivatives, with holomorphic coefficients, and for which one can define the notion of degree, as previously. Then if the two differential equations in question are equivalent under point transformations of the plane, then it is necessary that the two differential equations have the same degree, and that the corresponding curves on the cone $\mathscr S_3$ are projectively equivalent.
\end{proposition}
\begin{proof}
This follows because if two such differential equations are equivalent under point transformations, their first and second derivatives are related by an identity of the form \eqref{pro3}. Then equation \eqref{pro8} enables us to conclude.
\end{proof}


\end{document}